\documentclass[11pt,reqno]{amsart}
\usepackage[T1]{fontenc}
\usepackage[utf8]{inputenc}
\usepackage{lmodern}
\usepackage{amsmath,amssymb,amsthm,mathtools}
\usepackage[margin=1in]{geometry}
\usepackage[pdftex,pagebackref=true,backref=true,bookmarksopen=true]{hyperref}
\usepackage{float}\usepackage{placeins}\usepackage{graphicx}\usepackage{tikz}\usetikzlibrary{arrows.meta,patterns}
\usepackage[en-US]{datetime2} \DTMsettimestyle{iso} \date{\today}

\newtheorem{theorem}{Theorem}

\newtheorem{lemma}[theorem]{Lemma}
\newcommand{\R}{\mathbb R}
\newcommand{\Sn}{\mathbb S^{d-1}}

\begin{document}

\title{Radial conditions for global compact finite propagation for divergence-form elliptic operators}
\author{M.A.~Perelmuter}
\address{SCAD Soft Ltd., 3a Osvity Street, Kyiv, 03037, Ukraine}
\email{mikeperelmuter@gmail.com}

\begin{abstract}
We study global compact finite propagation (GCFP) for the wave equation $u_{tt}=\operatorname{div}(A(x)\nabla u)$ on $\mathbb R^d$, that is, the property that a disturbance starting in a compact set cannot reach spatial infinity in finite time. For continuous strictly positive matrices $A$ we show that GCFP implies $\int\limits_1^\infty a_*(r)^{-1/2}\,\mathrm{d}r=\infty$, where $a_*(r)=\min_{|x|=r}\langle A(x)x/|x|,x/|x|\rangle$. The proof uses heat-semigroup asymptotics. We also prove that the analogous divergence condition with $a^*(r)$, the maximum over the sphere, is sufficient even for measurable coefficients $A\in L^1_{\mathrm{loc}}$ with $A\ge I$. Finally, we construct a locally bounded measurable matrix with $a_*=a^*$ for which GCFP holds although the integral converges, so the necessary condition is not stable under the loss of continuity.
\end{abstract}

\maketitle

\section{Introduction}\label{sec:introduction}
Consider the Cauchy problem for the divergence-form wave equation
$$u_{tt}=\operatorname{div}(A(x)\nabla u),\qquad u(0)=\phi,\quad u_t(0)=\psi,$$
on $\R^d$, with compactly supported initial data. For smooth uniformly elliptic matrices with bounded coefficients, the classical finite propagation property holds with a uniform bound on the propagation speed\footnote{In the literature originally published in Russian~\cite{O2,O,RB,PS} the term \emph{finite propagation rate} is used for closely related properties.}. When the coefficients may grow at infinity, such a uniform bound need not exist. Two notions of finite propagation should therefore be distinguished. A \emph{uniform finite speed} means that supports expand at most at a fixed Euclidean speed, whereas \emph{global compact finite propagation} (GCFP) means that a disturbance starting in a fixed compact set cannot reach spatial infinity in finite time.

Finite propagation is closely connected with spectral theory and the analysis of functions of elliptic operators. The wave-equation method underlies kernel estimates on complete Riemannian manifolds~\cite{CGT} and is related to Gaussian bounds and Riesz transform~\cite{Sik}. Related notions of finite and infinite wave propagation speed on metric measure spaces, with applications to fractals, were studied by Ngai, Tang, and Xie~\cite{NTX}. Finite-propagation methods have also played an important role in proving essential self-adjointness of differential operators (see Povzner~\cite{Pov}, Berezanskii~\cite[\S 2.6]{Ber}, Chernoff~\cite{C,C2}, Kato~\cite{Kato}, Rauch--Taylor~\cite{RT}, Chumak~\cite{Ch}, Orochko~\cite{O2,O}, and Perelmuter~\cite{P}).

A substantial literature is devoted to sufficient conditions for finite propagation, including the work of Orochko~\cite{O2}, Perelmuter--Semenov~\cite{PS}, and Rofe-Beketov~\cite{RB}, as well as related developments cited above. Necessary conditions for global compact finite propagation have received less attention. A particularly relevant result is due to Rofe-Beketov~\cite[Theorem~1]{RB}, who established, for smooth strictly positive coefficients, an equivalence between global finite propagation and completeness of the Riemannian metric associated with $A^{-1}$.

The main results of this note concern \emph{radial} necessary and sufficient conditions for GCFP. For continuous strictly positive matrices, we show that GCFP implies
$$\int\limits_1^\infty\frac{\mathrm{d}r}{\sqrt{a_*(r)}}=\infty,\qquad a_*(r)=\min_{|\omega|=1}\left\langle A(r\omega)\omega,\omega\right\rangle.$$
Conversely, a sufficient condition is
$$\int\limits_1^\infty\frac{\mathrm{d}r}{\sqrt{a^*(r)}}=\infty,\qquad a^*(r)=\operatorname*{ess\,sup}_{|\omega|=1}\left\langle A(r\omega)\omega,\omega\right\rangle.$$
The quantities $\sqrt{a_*(r)}$ and $\sqrt{a^*(r)}$ represent the minimum and maximum local normal propagation speeds across the sphere of radius $r$. For smooth coefficients the sufficient condition is contained in the work of Rofe-Beketov~\cite[Theorem~2 and Example~1]{RB}. The sufficient condition does not require smoothness: in Section~\ref{sec:sufficient} we prove it for measurable coefficients $A\in L^1_{\mathrm{loc}}\left(\R^d;\operatorname{Sym}_d\right)$ satisfying $A\ge I$ almost everywhere, with the wave evolution defined by the closure of the quadratic form on $C_c^\infty\left(\R^d\right)$. Sufficient conditions for finite propagation speed in the case of nonsmooth coefficients were also investigated by Perelmuter--Semenov~\cite{PS}.

Our main objective is to investigate the necessity of the radial condition for GCFP. We prove it directly for continuous coefficients in Section~\ref{sec:necessary}, without the smoothness assumption of Rofe-Beketov's theorem. We then give a locally bounded measurable positive-definite matrix for which GCFP holds although the integral involving $a_*(r)$ converges. Thus the necessary condition is not stable under the loss of continuity.

The paper is organized as follows. Section~\ref{sec:setting} fixes the setting in both the continuous and the measurable cases. Section~\ref{sec:sufficient} contains the sufficient condition. Section~\ref{sec:necessary} establishes the necessary condition for continuous coefficients using heat-semigroup asymptotics. We then present the measurable counterexample. This failure of necessity came as a surprise to the author and motivated the present note.

\section{Setting and meaning of global finite propagation}\label{sec:setting}

\paragraph{\emph{Continuous coefficients.}}
Let $d\geq2$ and let $A\in C\left(\R^d;\operatorname{Sym}_d\right)$ be strictly positive definite at every point; no uniform upper or lower bounds at infinity are assumed. Consider the wave equation
$$u_{tt}=\operatorname{div}(A(x)\nabla u).$$
Write $g_x(v,w)=\langle A(x)^{-1}v,w\rangle$, $v,w\in\R^d$. For an absolutely continuous curve $\gamma:[a,b]\to\R^d$, define its $g$-length by
$$ L_g(\gamma)=\int\limits_a^b \langle A(\gamma(t))^{-1}\gamma'(t),\gamma'(t)\rangle^{1/2}\,\mathrm dt.$$
Let $d_g(x,y)$ be the infimum of these lengths over absolutely
continuous curves joining $x$ to $y$. Since $A$ is continuous and positive definite, the metric $d_g$ is locally equivalent to the Euclidean metric. Hence the $d_g$-topology agrees with the Euclidean topology. Consider the form
$$ \mathfrak a_0[u,v]=\int\limits_{\R^d} \langle A\nabla u,\nabla\overline v\rangle\,\mathrm dx, \qquad u,v\in C_c^\infty(\R^d).$$
Since $\lambda_K I\le A\le\Lambda_K I$ on each compact set $K$ for some $0<\lambda_K\le\Lambda_K<\infty$, the closability argument for measurable coefficients given below applies locally. Hence $\mathfrak a_0$ is closable. Let $\mathfrak a$ be its closure and let $H$ be the nonnegative self-adjoint operator associated with $\mathfrak a$. Set $W(t)=\cos(t\sqrt H)$ and $S(t)=s_t(H)$, where $s_t(\lambda)=\sin(t\sqrt\lambda)/\sqrt\lambda$ for $\lambda>0$ and $s_t(0)=t$. We say that $W$ has
\emph{global compact finite propagation} (GCFP) if, for every compact $K\subset\R^d$ and $T>0$, there is a compact $K_T$ such that
$$ \operatorname{supp}W(t)u\subset K_T\qquad (|t|\le T)$$
for all $u\in L^2$ supported in $K$. The analogous property for $S(t)$ follows from $S(t)=\int\limits_0^t W(s)\,\mathrm ds$. Throughout, the wave evolution is associated with the operator of the closed form $\mathfrak a$.

\paragraph{\emph{Measurable coefficients.}}
Now assume that $A:\R^d\to\operatorname{Sym}_d$ is measurable, finite almost everywhere, and
$$A\in L^1_{\mathrm{loc}}\left(\R^d\right),\qquad A(x)\ge I\quad\text{a.e.}$$
The form $\mathfrak a_0[u,v]=\int\limits\left\langle A\nabla u,\nabla\overline v\right\rangle\,\mathrm{d}x$ on $C_c^\infty\left(\R^d\right)$ is densely defined and closable. Indeed, if $u_n\to0$ in $L^2$ and $A^{1/2}\nabla u_n$ is Cauchy in $L^2$, then $\nabla u_n$ is Cauchy in $L^2$ because $A\ge I$, and its limit is zero. Passing to a subsequence along which both $\nabla u_n$ and $A^{1/2}\nabla u_n$ converge almost everywhere, we conclude that the limit of $A^{1/2}\nabla u_n$ is zero as well.
Let $\mathfrak a$ denote the closure of $\mathfrak a_0$, and let $H$ be the nonnegative self-adjoint operator associated with $\mathfrak a$. As in the continuous setting, $H$ is defined as the operator of the closed form. Set $W(t)=\cos\left(t\sqrt H\right)$ and $S(t)=s_t(H)$, where $s_t(\lambda)=\sin\left(t\sqrt\lambda\right)/\sqrt\lambda$ for $\lambda>0$ and $s_t(0)=t$. Thus $S(t)$ is bounded for each $t$. We interpret GCFP for this realization as follows: for every compact set $K$ and every $T>0$, there exists a compact set $K_T$ such that $W(t)f$ and $S(t)g$ are supported in $K_T$ whenever $f$ and $g$ are supported in $K$ and $|t|\le T$.
In both settings the form $\mathfrak a$ is a Dirichlet form. Indeed, for real-valued $u\in C_c^\infty(\R^d)$ the function $v=\min(1,\max(0,u))$ is compactly supported and Lipschitz. Hence it belongs to $D(\mathfrak a)$ (see the argument at the beginning of the proof of Lemma~\ref{lem:form-multiplier}). The Sobolev chain rule gives $\langle A\nabla v,\nabla v\rangle\le\langle A\nabla u,\nabla u\rangle$ a.e., so $\mathfrak a[v]\le\mathfrak a[u]$.
For real-valued $u\in D(\mathfrak a)$, take real-valued $u_n\in C_c^\infty(\R^d)$ with $u_n\to u$ in the form norm. Then $\min(1,\max(0,u_n))$ is bounded in the form norm and converges to $\min(1,\max(0,u))$ in $L^2$. By weak compactness and lower semicontinuity, the latter belongs to $D(\mathfrak a)$ and $\mathfrak a[\min(1,\max(0,u))]\le\mathfrak a[u]$. Consequently, $e^{-sH}$ is positivity preserving and the product and truncation properties of~\cite[Theorem~1.4.2]{FOT} are available.

\paragraph{\emph{Radial quantities.}}
For $r>0$ and $\omega\in\Sn$ put $a(r,\omega)=\langle A(r\omega)\omega,\omega\rangle$. For continuous $A$ define
$$a_*(r)=\min_{\omega\in\Sn}a(r,\omega),\qquad a^*(r)=\max_{\omega\in\Sn}a(r,\omega).$$
Both functions are then positive and continuous. For measurable $A$ define, for almost every $r>0$,
$$a_*(r)=\operatorname*{ess\,inf}_{\omega\in\Sn}\langle A(r\omega)\omega,\omega\rangle,\qquad a^*(r)=\operatorname*{ess\,sup}_{\omega\in\Sn}\langle A(r\omega)\omega,\omega\rangle,$$
with the convention $1/\sqrt{+\infty}=0$. For continuous $A$ the pointwise and essential definitions agree.

The function $a(r,\omega)$ is jointly measurable. The polar-coordinate formula
$$\int\limits_{\R^d}f(x)\,\mathrm{d}x=\int\limits_0^\infty\int\limits_{\Sn}f(r\omega)r^{d-1}\,\mathrm{d}\sigma(\omega)\,\mathrm{d}r,$$
valid for every nonnegative measurable function $f$, and Tonelli's theorem~\cite[Theorem~2.37]{Folland} imply that, for almost every $r>0$, the angular section $a(r,\cdot)$ is finite almost everywhere and satisfies $a(r,\omega)\ge 1$ almost everywhere on $\Sn$. Moreover, for every $t\in\R$, the functions
$r\mapsto \sigma\left\{\omega\in \Sn:a(r,\omega)>t\right\}$ and $r\mapsto \sigma\left\{\omega\in \Sn:a(r,\omega)<t\right\}$ are measurable by Tonelli's theorem. Since the conditions $a^*(r)\le t$ and $a_*(r)\ge t$ are equivalent to the vanishing of the respective level-set measures, both essential bounds are measurable. As noted above, $a(r,\cdot)$ is finite almost everywhere for almost every $r>0$. Hence $a_*$ is finite almost everywhere, whereas $a^*$ may take the value $+\infty$.

\section{Sufficient condition}\label{sec:sufficient}
We begin with the following abstract lemma.
\begin{lemma}\label{lem:abstract}
Let $H\ge0$ be a self-adjoint operator in a Hilbert space $\mathcal H$, let $f,g\in\mathcal H$ and $\rho>0$. Suppose that
$$ \left|\left\langle e^{-tH}f,g\right\rangle\right|\le e^{-\rho^2/4t}\|f\|\,\|g\|\qquad\text{for all }t>0.$$
Then $\langle\cos\left(\tau\sqrt H\right)f,g\rangle=0$ for $|\tau|<\rho$.
\end{lemma}
\begin{proof}
The function $F(z)=\left\langle e^{-zH}f,g\right\rangle$ is analytic in the right half-plane $\operatorname{Re}z>0$, extends continuously to its boundary, and satisfies $|F(z)|\le M:=\|f\|\|g\|$ for $\operatorname{Re}z\ge0$. By hypothesis, $|F(t)|\le M e^{-\rho^2/(4t)}$ for $t>0$. Applying the Phragm\'en--Lindel\"of principle (see \cite[Lemma~4.2]{SW}) for sectors of opening angle $\pi/2$ to $F(1/w)e^{\rho^2w/4}$ in the upper and lower quadrants of the right half-plane (the function is bounded near the vertex, bounded by $M$ on the boundary rays, and has growth of order at most one) yields
$$|F(z)|\le M\exp\left(-\frac{\rho^2}{4}\operatorname{Re}\frac1z\right),\qquad \operatorname{Re}z>0.$$
The passage from this bound to $\left\langle\cos\left(\tau\sqrt H\right)f,g\right\rangle=0$ for $|\tau|<\rho$ follows the Laplace-transform argument in the proof of~\cite[Theorem~3.4]{CS}, which applies in the present scalar setting.
\end{proof}
We denote by $D(\mathfrak a)_{\mathrm{loc}}$ the local form domain, consisting of all $u\in L^2_{\mathrm{loc}}\left(\R^d\right)$ such that, for every compact set $K\subset\R^d$, there exists $v\in D(\mathfrak a)$ with $u=v$ almost everywhere on $K$.
\begin{lemma}\label{lem:form-multiplier}
Let $(\mathfrak a,D(\mathfrak a))$ be the closed form of Section~\ref{sec:setting}, either for measurable $A\in L^1_{\mathrm{loc}}$ with $A\ge I$ almost everywhere or for continuous strictly positive definite $A$. Suppose $\Phi:\R^d\to[0,\infty)$ is continuous and satisfies $\Phi(x)\to\infty$ as $|x|\to\infty$, and every truncation $\Phi_N=\min(\Phi,N)$ has the following properties: $\Phi_N\in D(\mathfrak a)_{\mathrm{loc}}$, and $\left\langle A\nabla\Phi_N,\nabla\Phi_N\right\rangle\le C^2$ almost everywhere, with $C$ independent of $N$. Then the associated cosine family $W(t)$ has GCFP.
\end{lemma}
\begin{proof}
If $u$ is compactly supported and Lipschitz, its standard mollifications $u_\varepsilon=u*\eta_\varepsilon$, where $\eta_\varepsilon\in C_c^\infty(\R^d)$, belong to $C_c^\infty(\R^d)$ and converge to $u$ in $L^2$, while $\nabla u_\varepsilon\to\nabla u$ almost everywhere and $|\nabla u_\varepsilon|\le\|\nabla u\|_\infty$. Since the gradients have a common compact support for sufficiently small $\varepsilon$, the dominated convergence theorem and the assumption $A\in L^1_{\mathrm{loc}}$ imply convergence in the form norm. This argument also applies to continuous strictly positive definite $A$, which is bounded above and below on each compact set. Hence $u\in D(\mathfrak a)$. For $\alpha\in\R$, put $w=e^{\pm\alpha\Phi_N}$. The function $w$ is bounded and belongs to $D(\mathfrak a)_{\mathrm{loc}}$; the Sobolev chain rule and the assumed bound on $\langle A\nabla\Phi_N,\nabla\Phi_N\rangle$ give
$$\langle A\nabla w,\nabla w\rangle\le \alpha^2 e^{2|\alpha|N}C^2\quad\text{a.e.}$$
For $u\in C_c^\infty$, the local form-domain property of $w$, together with the truncation and product properties of Dirichlet forms (see~\cite[Theorem~1.4.2(ii), (iii)]{FOT}), implies that $wu\in D(\mathfrak a)$. Indeed, there exists a bounded function $v\in D(\mathfrak a)$ that coincides with $w$ almost everywhere in a neighborhood of $\operatorname{supp}u$. The product rule for weak derivatives then applies locally. Moreover,
$$\mathfrak a[wu]\le 2e^{2|\alpha|N}\mathfrak a[u]+2\alpha^2e^{2|\alpha|N}C^2\|u\|_2^2.$$
Thus multiplication by $w$ extends continuously from the form core to $D(\mathfrak a)$, and the product rule extends by form-norm approximation. Since $A$ is real symmetric and $\Phi_N$ is real, the mixed terms are purely imaginary and disappear on taking real parts, which gives
$$\operatorname{Re}\mathfrak a\left(e^{-\alpha\Phi_N}u,e^{\alpha\Phi_N}u\right)
=\mathfrak a[u]-\alpha^2\int\limits |u|^2\langle A\nabla\Phi_N,\nabla\Phi_N\rangle\,\mathrm{d}x
\ge-\alpha^2C^2\|u\|_2^2.$$
 The Davies exponential perturbation method (see Davies~\cite[Section~3.2]{Davies}) yields the Davies--Gaffney estimate (cf.~\cite[(3.2)]{CS}). Applying it with the admissible weight $\Phi_N/C$ and $N>\max\{|\ell_1|,|\ell_2|\}$ gives, for $f,g\in L^2\left(\R^d\right)$ supported in $\{\Phi\le \ell_1\}$ and $\{\Phi\ge \ell_2\}$, respectively, with $\ell_1<\ell_2$,
$$\left|\left\langle e^{-tH}f,g\right\rangle\right|\le \exp\!\left(-\frac{(\ell_2-\ell_1)^2}{4C^2t}\right)\|f\|_2\|g\|_2,\qquad t>0.$$ Applying Lemma~\ref{lem:abstract} to this fixed pair $f,g$ with $\rho=(\ell_2-\ell_1)/C$, we obtain
$$ \left\langle W(t)f,g\right\rangle=0\qquad\text{if }|t|<(\ell_2-\ell_1)/C.$$
The same support restriction holds for $S(t)=\int\limits_0^tW(s)\,\mathrm{d}s$. Since $\Phi(x)\to\infty$ as $|x|\to\infty$, every sublevel set $\{x:\Phi(x)\le M\}$ is compact. Thus for initial data supported in a fixed compact set $K$ and $|t|\le T$, both evolutions have support in the compact set $\{\Phi\le \max_K\Phi+CT\}$, proving GCFP.
\end{proof}

\begin{theorem}\label{thm:sufficient}
Let $A$ be as in the measurable setting of Section~\ref{sec:setting} (in particular, $A\in L^1_{\mathrm{loc}}$ and $A\ge I$ a.e.). If
$$\int\limits_1^\infty\frac{\mathrm{d}r}{\sqrt{a^*(r)}}=\infty,$$
then the wave evolution has GCFP.
\end{theorem}
\begin{proof}
Put $\Phi(x)=0$ for $|x|\le1$ and $\Phi(x)=\int\limits_1^{|x|}a^*(s)^{-1/2}\,ds$ for $|x|>1$. Since $a^*\ge1$, we have $0\le a^*(s)^{-1/2}\le1$, so $\Phi$ is continuous and Lipschitz, and so are its truncations $\Phi_N=\min(\Phi,N)$. Since $A\in L^1_{\mathrm{loc}}$, the truncations belong locally to the closed-form domain (see the first paragraph of the proof of Lemma~\ref{lem:form-multiplier}). Almost everywhere, $\nabla\Phi=a^*(|x|)^{-1/2}\,x/|x|$, hence
$$\langle A\nabla\Phi,\nabla\Phi\rangle=\frac{a(r,\omega)}{a^*(r)}\le1$$
(with the convention $1/\sqrt{+\infty}=0$ on radii with $a^*=\infty$, so that $\nabla\Phi=0$ there). The same bound holds for $\Phi_N$. Divergence of the integral ensures that $\Phi(x)\to\infty$ as $|x|\to\infty$, so Lemma~\ref{lem:form-multiplier} yields GCFP.
\end{proof}

\section{Necessary condition}\label{sec:necessary}

\subsection{Continuous coefficients}\label{sec:continuous}

Let $A\in C(\R^d;\operatorname{Sym}_d)$ be strictly positive definite at every point. Let $\mathfrak a$ be the closure of $\int\limits_{\R^d}\langle A\nabla u,\nabla\bar u\rangle\,\mathrm dx$ on $C_c^\infty(\R^d)$, let $H$ be its associated self-adjoint operator, and put $W(t)=\cos(t\sqrt H)$. We use the length distance $d_g$ of
$g=A^{-1}$ defined in Section~\ref{sec:setting}, and write $d_g(U,V)=\inf_{x\in U,y\in V}d_g(x,y)$.

\begin{lemma}\label{lem:propagation-distance}
Assume that $W$ has GCFP. Let $U\subset\R^d$ be a nonempty relatively compact open set, let $T>0$, and let $K_T$ be a compact set containing $\operatorname{supp}W(t)f$ for all $f$ supported in $\overline U$ and
$|t|\le T$. Then $T\le d_g(U,V)$ for every nonempty relatively compact open set $V$ with $V\cap K_T=\emptyset$.
\end{lemma}
\begin{proof}
For such $V$ we have $\langle W(t)\mathbf1_U,\mathbf1_V\rangle=0$ for $|t|\le T$. By the Gaussian representation of the heat semigroup,
$$ e^{-sH}=\frac1{\sqrt{4\pi s}}\int\limits_{\R} e^{-t^2/(4s)}W(t)\,\mathrm dt$$
and $\|W(t)\|_{2\to2}\le1$ give
$$ 0\le\langle e^{-sH}\mathbf1_U,\mathbf1_V\rangle \le\frac{|U|^{1/2}|V|^{1/2}}{\sqrt{4\pi s}}
 \int\limits_{|t|>T}e^{-t^2/(4s)}\,\mathrm dt.$$
Thus, with $\log0=-\infty$,
\begin{equation}\label{eq:continuous-wave-bound}
 \limsup_{s\downarrow0}4s\log
 \langle e^{-sH}\mathbf1_U,\mathbf1_V\rangle\le-T^2.
\end{equation}

The form $\mathfrak a$ is strongly local, since
$$\mathfrak a[u,v]=\int\limits_{\R^d}\langle A\nabla u,\nabla\overline v\rangle\,\mathrm dx,\qquad u,v\in D(\mathfrak a).$$

Every $u\in D(\mathfrak a)_{\mathrm{loc}}$ belongs to $H^1_{\mathrm{loc}}$ by local ellipticity. For bounded real-valued such $u$ and $h\in C_c^\infty$, the localized energy functional of~\cite{AH} is
$$ I_u(h)=2\mathfrak a(uh,u)-\mathfrak a(u^2,h)
 =2\int\limits_{\R^d}h\langle A\nabla u,\nabla u\rangle\,\mathrm dx.$$
Here the form expressions are interpreted locally on a neighborhood of $\operatorname{supp}h$. Hence every $u$ in the class $\mathcal D_0$ (see~\cite[Definition~2.6]{AH}) satisfies $\langle A\nabla u,\nabla u\rangle\le1/2$ almost everywhere, by testing with nonnegative $h\in C_c^\infty$. Local ellipticity gives
$\nabla u\in L^\infty_{\mathrm{loc}}$. Hence $u$ coincides almost everywhere with a locally Lipschitz function.

Denoting this locally Lipschitz function again by $u$, we have $|u(x)-u(y)|\le d_g(x,y)/\sqrt2$. Indeed, let $\gamma$ be an absolutely continuous curve from $x$ to $y$. On a compact neighborhood of the curve $\gamma$, uniform continuity and local ellipticity of $A$ imply $A(z)\le(1+\delta)^2 A(w)$ whenever $|z-w|$ is sufficiently small. By convexity of
$\xi\mapsto\langle A(z)\xi,\xi\rangle^{1/2}$, the standard mollifications $u_\varepsilon$ satisfy
$$
\langle A(z)\nabla u_\varepsilon(z),
\nabla u_\varepsilon(z)\rangle^{1/2}
\le (1+\delta)/\sqrt2
$$
at every point $z$ of the curve $\gamma$ for
sufficiently small $\varepsilon$. By the chain rule and the Cauchy--Schwarz inequality, taking $\gamma:[0,1]\to\R^d$, we obtain
\begin{align*}
|u_\varepsilon(y)-u_\varepsilon(x)|=\left|\int\limits_0^1\frac{\mathrm d}{\mathrm dt}  u_\varepsilon(\gamma(t))\,\mathrm dt\right|
\le\int\limits_0^1|\langle\nabla u_\varepsilon(\gamma(t)),\gamma'(t)\rangle|\,\mathrm dt\\
\le\frac{1+\delta}{\sqrt2}\int\limits_0^1
\langle A(\gamma(t))^{-1}\gamma'(t),\gamma'(t)\rangle^{1/2}\,\mathrm dt =\frac{1+\delta}{\sqrt2}L_g(\gamma).
\end{align*}
Letting $\varepsilon\to0$ and then $\delta\to0$,
and taking the infimum over all such curves, we obtain
$$ |u(x)-u(y)|\le d_g(x,y)/\sqrt2.$$

The distance of Ariyoshi--Hino is defined by
$$d_{\rm AH}(U,V)= \sup_{u\in\mathcal D_0} \left( \operatorname*{ess\,inf}_{y\in V}u(y) -\operatorname*{ess\,sup}_{x\in U}u(x)\right)$$
(see \cite[p.~1241]{AH}).
Since $U$ and $V$ are open and $u$ is continuous,
the essential infimum and supremum coincide with
the ordinary infimum and supremum, respectively.
Hence
$$\begin{aligned}
\operatorname*{ess\,inf}_{y\in V}u(y)
-\operatorname*{ess\,sup}_{x\in U}u(x)
=\inf_{x\in U,\,y\in V}(u(y)-u(x))\\
\le\frac{1}{\sqrt2}\inf_{x\in U,\,y\in V}d_g(x,y)=\frac{d_g(U,V)}{\sqrt2}.
\end{aligned}$$
Taking the supremum over $u\in\mathcal D_0$, we obtain
$$d_{\rm AH}(U,V)\le\frac{d_g(U,V)}{\sqrt2}.$$
The theorem of Ariyoshi--Hino \cite[Theorem~2.7]{AH} gives
$$ \lim_{s\downarrow0}4s\log\langle e^{-sH}\mathbf1_U,\mathbf1_V\rangle
 =-2d_{\rm AH}(U,V)^2\ge-d_g(U,V)^2.$$
Together with~\eqref{eq:continuous-wave-bound}, this proves $T^2\le d_g(U,V)^2$.
\end{proof}

\begin{theorem}\label{thm:continuous-necessary}
Let $A\in C(\R^d;\operatorname{Sym}_d)$ be strictly positive definite. If the associated cosine family has GCFP, then
$$ \int\limits_1^\infty\frac{\mathrm dr}{\sqrt{a_*(r)}}=\infty, \qquad a_*(r)=\min_{\omega\in\Sn} \langle A(r\omega)\omega,\omega\rangle.$$
\end{theorem}
\begin{proof}
Suppose the integral is finite. Let
$$ L=\int\limits_1^\infty\frac{\mathrm dr}{\sqrt{a_*(r)}}<\infty.$$
Decompose $A(r\omega)\omega=a(r,\omega)\omega+b(r,\omega)$, where $b(r,\omega)\perp\omega$. For $r\ge1$ and
$z\in\R^d\setminus\{0\}$ define
$$ F(r,z)=\frac{b(r,z/|z|)}{r\,a(r,z/|z|)}.$$
This vector field is continuous and tangent to the spheres $|z|=\mathrm{const}$. It is bounded on each set
$[1,R]\times\{z:1/2\le|z|\le2\}$. By Peano's theorem (see \cite[Chapter~1, Theorem~1.2]{CL}), there exists a local solution
of $\omega'=F(r,\omega)$ with $|\omega(1)|=1$. Since $z\cdot F(r,z)=0$, we have $|\omega(r)|=1$ throughout its interval of existence.
The boundedness of $F$ on the above sets and the continuation theorem (see \cite[Chapters~1--2]{CL}) imply that a maximal solution is defined for all $r\ge1$. Thus we obtain a solution $\omega:[1,\infty)\to\Sn$.

Set $x(r)=r\omega(r)$. Then
$$ x'(r)=\frac{A(x(r))\omega(r)}{a(r,\omega(r))},  \qquad  \langle A(x(r))^{-1}x'(r),x'(r)\rangle  =\frac1{a(r,\omega(r))}  \le\frac1{a_*(r)}.$$
Thus $|x(r)|=r$, and the $g$-length of $x$ on $[1,\infty)$ is at most $L$.

Assume GCFP. Choose a relatively compact open ball $U$ containing $x(1)$, let $T=L+1$, and let $K_T$ be given by
Lemma~\ref{lem:propagation-distance}. Choose $R>\max_{z\in K_T}|z|$ and a relatively compact open ball
$V$ containing $x(R)$ and disjoint from $K_T$. Then
$$ d_g(U,V)\le d_g(x(1),x(R))\le L<T,$$
contrary to Lemma~\ref{lem:propagation-distance}.
\end{proof}

\subsection{Failure of radial necessity for measurable coefficients}
The radial condition need not be necessary for GCFP of the wave evolution associated with the corresponding Dirichlet form when the coefficients are merely measurable. We exhibit this failure below, in the measurable setting of Section~\ref{sec:setting}. The idea is to construct a coefficient matrix with a rapidly growing radial component but sufficiently strong anisotropy to ensure GCFP despite the failure of the radial condition.

\begin{theorem}\label{thm:rough-counterexample}
There exists $A\in L^\infty_{\mathrm{loc}}\left(\R^2;\operatorname{Sym}_2\right)$ with $A\ge I$ such that, with respect to the orthonormal basis $(e_r,e_\theta)$, where $e_r$ and $e_\theta$ are the radial and angular unit vectors, for $r>0$,
$$ A(r,\theta)=\begin{pmatrix}
 1+r^4 & -r^5 h'(\theta)\\
 -r^5 h'(\theta) & 1+r^6
 \end{pmatrix},\qquad
 h(\theta)=\operatorname{dist}_{\mathbb S^1}(\theta,0), $$
where $h'$ is defined to be zero at the angles $\theta=0,\pi$ (modulo $2\pi$), at which its classical derivative does not exist. The definition of $A$ is understood almost everywhere. Moreover,
$$a_*(r)=a^*(r)=1+r^4\qquad(r>0),$$
while the associated cosine family $W(t)$ has GCFP. Consequently, the necessity of $\int\limits_1^\infty a_*(r)^{-1/2}\,\mathrm{d}r=\infty$ does not extend to arbitrary measurable coefficients, even when $a_*(r)=a^*(r)$.
\end{theorem}
\begin{proof}
The function $h$ is $2\pi$-periodic and Lipschitz, with $|h'|=1$ almost everywhere.\footnote{No $C^1$ periodic function $h$ can satisfy $|h'|=1$ almost everywhere, since its derivative must vanish at a point of maximum or minimum. Thus, the lack of continuous differentiability is essential for the present construction.} For $r>0$ put
$$v(r,\theta)=e_r-rh'(\theta)e_\theta,\qquad A(r,\theta)=I+r^4v(r,\theta)\otimes v(r,\theta).$$
Since $|v|^2=1+r^2$, we have $I\le A\le\bigl(1+r^4(1+r^2)\bigr)I$, so $A\ge I$ and $A$ is locally essentially bounded. Since $\langle v,e_r\rangle=1$, we have $\langle Ae_r,e_r\rangle=1+r^4$ for all $r>0$, and the integral of its reciprocal square root is finite.

Consider $f(r,\theta)=r+h(\theta)$ for $r>0$ (the function $f$ provides a distance-like coordinate in which propagation towards infinity remains uniformly controlled, despite the rapid growth of the radial coefficient). Almost everywhere,
$$\nabla f=e_r+\frac{h'(\theta)}r e_\theta,\qquad \langle v,\nabla f\rangle=1-(h')^2=0,$$
whence
$$\langle A\nabla f,\nabla f\rangle=|\nabla f|^2=1+r^{-2}.$$
This expression is bounded for $r\ge2$ but tends to infinity as $r\to0$. We therefore cut $f$ off near the origin. Fix $\eta\in C^\infty([0,\infty))$ with $0\le\eta\le1$, $\eta=0$ on $[0,1]$ and $\eta=1$ on $[2,\infty)$, and put $F(x)=\eta(|x|)f(x)$. Then $F\ge0$ is locally Lipschitz on $\R^2$ and $F=f$ for $|x|\ge2$. For $|x|\ge2$ we have $\nabla F=\nabla f$ almost everywhere, since $\eta=1$ and $\eta'=0$. Hence
$$\langle A\nabla F,\nabla F\rangle=1+r^{-2}\le5/4.$$
For $1\le|x|\le2$ the gradient $\nabla F$ is bounded and $A$ is locally essentially bounded. For $|x|\le1$ we have $F=0$. Hence $\langle A\nabla F,\nabla F\rangle\le C<\infty$ almost everywhere, and the same bound holds for the truncations $F_N=\min(F,N)$. Since $A$ is locally essentially bounded, $F$ and its truncations belong to the local domain of the closed form $\mathfrak a$. Moreover $F(x)=|x|+h(\theta)\ge|x|$ for $|x|\ge2$, so $F(x)\to\infty$ as $|x|\to\infty$. Lemma~\ref{lem:form-multiplier}, applied to $\Phi=F$, therefore ensures that the supports of the waves remain in a fixed compact set for bounded times. This proves GCFP.
\end{proof}

A geometric illustration of why the matrix of Theorem~\ref{thm:rough-counterexample} does not push disturbances to infinity is given in Appendix~\ref{app:geometry}.

\appendix
\section{Geometric picture of the counterexample}\label{app:geometry}

This appendix explains why the matrix of Theorem~\ref{thm:rough-counterexample} does not push disturbances to infinity although $a(r,\omega)=1+r^4$ for all $r>0$ and all $\omega$.

\subsection{Local picture}
For a point $x$ let
$$E_x=\{y\in\R^2:\ \langle A(x)^{-1}y,y\rangle\le1\}$$
be the pointwise unit ball of the metric tensor $g=A^{-1}$ (defined almost everywhere); it is the set of displacements covered in unit time at $x$. Its support function is $\max_{y\in E_x}\langle\xi,y\rangle=\langle A(x)\xi,\xi\rangle^{1/2}$.

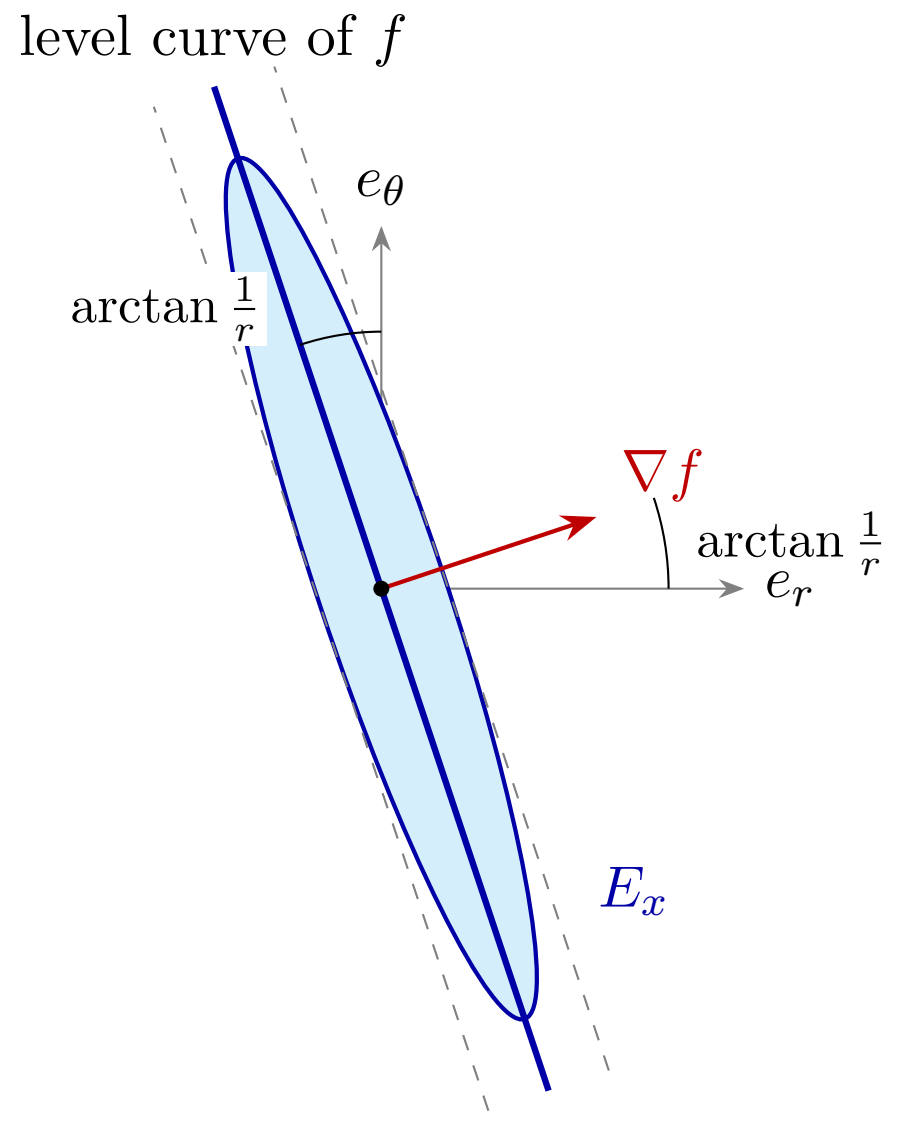
\begin{figure}[htbp]
\centering
\begin{tikzpicture}[>=Stealth]
 \draw[->,gray] (0,0)--(2.4,0) node[right,black] {$e_r$};
 \draw[->,gray] (0,0)--(0,2.4) node[above,black] {$e_\theta$};
 \begin{scope}[rotate=108.43]
  \draw[thick,blue!65!black,fill=cyan!15] (0,0) ellipse (3.0 and 0.42);
  \draw[very thick,blue!65!black] (-3.5,0)--(3.5,0);
  \draw[dashed,gray] (-3.5,0.42)--(3.5,0.42);
  \draw[dashed,gray] (-3.5,-0.42)--(3.5,-0.42);
  \draw[->,thick,red!75!black] (0,0)--(0,-1.5);
 \end{scope}
 \fill (0,0) circle (1.5pt);
 \node[right,red!75!black] at (1.45,0.75) {$\nabla f$};
 \node[above] at (-1.11,3.32) {level curve of $f$};
 \node[right,blue!65!black] at (1.3,-2.0) {$E_x$};
 \draw (0,1.7) arc (90:108.43:1.7);
 \node[left,fill=white,inner sep=1pt] at (-0.75,1.85) {\small $\arctan\frac1r$};
 \draw (1.9,0) arc (0:18.43:1.9);
 \node[right] at (1.95,0.3) {\small $\arctan\frac1r$};
\end{tikzpicture}
\caption{Indicatrix $E_x$ of the metric $A^{-1}$ at a point with $r=3$, $h'=1$, in the basis $(e_r,e_\theta)$. The long axis makes the angle $\arctan(1/r)$ with $e_\theta$, and $\nabla f$ makes the same angle with $e_r$.}
\label{fig:indicatrix}
\end{figure}

For the matrix of Theorem~\ref{thm:rough-counterexample} and $r>0$ the eigenvalues of $A=I+r^4v\otimes v$ are $1$ on $v^\perp$ and $1+r^4|v|^2=1+r^4(1+r^2)$ on $v$. Hence $E_x$ is a needle-shaped ellipse with semi-axes $\bigl(1+r^4(1+r^2)\bigr)^{1/2}\ge r^3$ along $v$ and $1$ across (Fig.~\ref{fig:indicatrix}). The key point is that $\nabla f$, where $f=r+h(\theta)$, is orthogonal to $v$, i.e.\ it points along the \emph{short} axis. Equivalently, $A\nabla f=\nabla f$, and therefore
$$\max_{y\in E_x}\langle\nabla f,y\rangle=\langle A\nabla f,\nabla f\rangle^{1/2}=|\nabla f|=\sqrt{1+r^{-2}}.$$
Here $\sqrt{1+r^{-2}}\le\sqrt{10/9}$ for $r\ge3$.
The long axis of $E_x$ is parallel to $v$, hence tangent to the level curve of $f$ through $x$.

\subsection{Global picture}
For $c\ge3+\pi$, the level curve
$$\Gamma_c=\{r+h(\theta)=c\}=\{r=c-|\theta|,\ |\theta|\le\pi\}$$
lies entirely in $r\ge3$. It consists of two Archimedean spiral arcs meeting at the corners $\theta=0$ and $\theta=\pi$, and is a closed piecewise smooth curve contained in $\{|x|\le c\}$. For smaller $c$, the level sets also depend on the cutoff of $f$ near the origin.

Away from the corners, the long axes of $E_x$ are tangent to $\Gamma_c$, since $\langle v,\nabla f\rangle=0$ there. For the Euclidean unit normal $n=\nabla f/|\nabla f|$, the maximal normal velocity is $\langle An,n\rangle^{1/2}=1$, whereas the maximal rate of change of $f$ under displacements in $E_x$ is $|\nabla f|\le\sqrt{10/9}$.

By Lemma~\ref{lem:form-multiplier} applied with $\Phi=F$, which coincides
with $f$ for $|x|\ge2$ (see the proof of Theorem~\ref{thm:rough-counterexample}),
data supported in a compact set $K$ produce a solution with
$$ \operatorname{supp}u(t)\subset \{\min_K F-C|t|\le F\le\max_K F+C|t|\}.$$
The lower bound follows by interchanging the two support sets in the Davies--Gaffney estimate, or equivalently by using the weight $-\Phi_N$. Thus, if $|F-F(p)|\le\varepsilon$ on the initial support, the disturbance
stays in the band $\{|F-F(p)|\le C|t|+\varepsilon\}$. Outside the disk $r<2$, for $|p|\ge2$, this is the band
$\{|f-f(p)|\le C|t|+\varepsilon\}$ around $\Gamma_{f(p)}$, although $E_x$ has length $\ge r^3$ along the level curves of $f$.

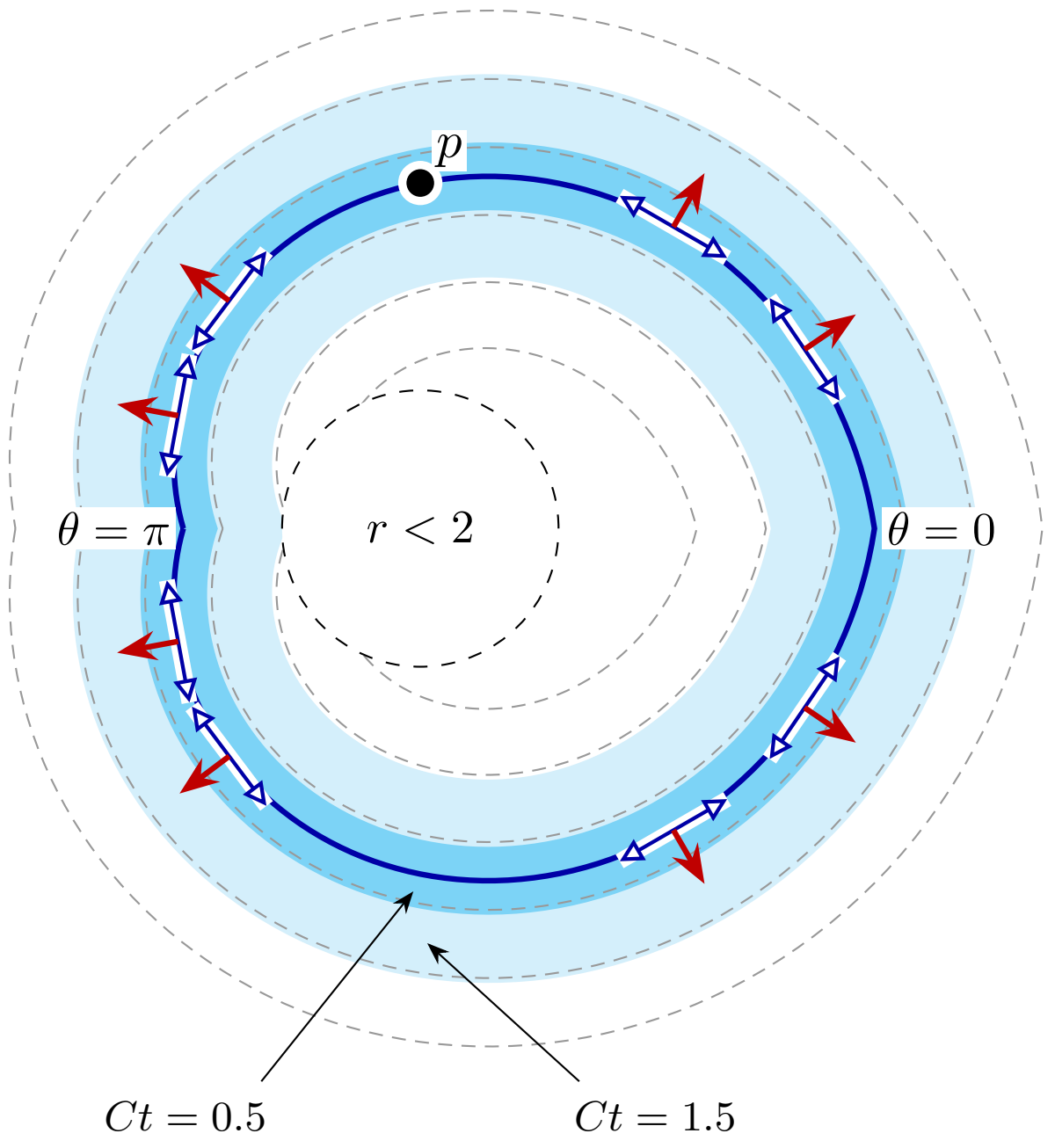
\begin{figure}[htbp]
\centering
\begin{tikzpicture}[scale=0.52,>=Stealth]
 \def\cc{6.5708}
 \fill[even odd rule,cyan!15] plot[domain=-180:180,samples=181,variable=\t]
   ({(\cc+1.5-abs(\t)*pi/180)*cos(\t)},{(\cc+1.5-abs(\t)*pi/180)*sin(\t)}) -- cycle
   plot[domain=-180:180,samples=181,variable=\t]
   ({(\cc-1.5-abs(\t)*pi/180)*cos(\t)},{(\cc-1.5-abs(\t)*pi/180)*sin(\t)}) -- cycle;
 \fill[even odd rule,cyan!45] plot[domain=-180:180,samples=181,variable=\t]
   ({(\cc+0.5-abs(\t)*pi/180)*cos(\t)},{(\cc+0.5-abs(\t)*pi/180)*sin(\t)}) -- cycle
   plot[domain=-180:180,samples=181,variable=\t]
   ({(\cc-0.5-abs(\t)*pi/180)*cos(\t)},{(\cc-0.5-abs(\t)*pi/180)*sin(\t)}) -- cycle;
 \foreach \c in {4,5,6,7,8,9}{
  \draw[gray!80,thin,densely dashed] plot[domain=-180:180,samples=181,variable=\t]
   ({(\c-abs(\t)*pi/180)*cos(\t)},{(\c-abs(\t)*pi/180)*sin(\t)});
 }
 \draw[very thick,blue!65!black] plot[domain=-180:180,samples=181,variable=\t]
   ({(\cc-abs(\t)*pi/180)*cos(\t)},{(\cc-abs(\t)*pi/180)*sin(\t)});
 \fill[white] (0,0) circle (2);
 \draw[dashed] (0,0) circle (2);
 \node at (0,0) {\small $r<2$};
 \foreach \t in {25,50,130,155,-25,-50,-130,-155}{
  \pgfmathsetmacro{\r}{\cc-abs(\t)*pi/180}
  \pgfmathsetmacro{\s}{sign(\t)}
  \pgfmathsetmacro{\px}{\r*cos(\t)}
  \pgfmathsetmacro{\py}{\r*sin(\t)}
  \pgfmathsetmacro{\tx}{-\s*cos(\t)-\r*sin(\t)}
  \pgfmathsetmacro{\ty}{-\s*sin(\t)+\r*cos(\t)}
  \pgfmathsetmacro{\n}{sqrt(\tx*\tx+\ty*\ty)}
  \pgfmathsetmacro{\gx}{cos(\t)-(\s/\r)*sin(\t)}
  \pgfmathsetmacro{\gy}{sin(\t)+(\s/\r)*cos(\t)}
  \pgfmathsetmacro{\gn}{sqrt(\gx*\gx+\gy*\gy)}
  \draw[white,line width=3.4pt] ({\px-0.9*\tx/\n},{\py-0.9*\ty/\n})--({\px+0.9*\tx/\n},{\py+0.9*\ty/\n});
  \draw[{Triangle[open,length=5pt,width=5pt]}-{Triangle[open,length=5pt,width=5pt]},line width=0.8pt,blue!65!black] ({\px-0.9*\tx/\n},{\py-0.9*\ty/\n})--({\px+0.9*\tx/\n},{\py+0.9*\ty/\n});
  \draw[->,line width=1.2pt,red!75!black] (\px,\py)--({\px+0.9*\gx/\gn},{\py+0.9*\gy/\gn});
 }
 \fill[white] (0,5) circle (0.32);
 \fill (0,5) circle (0.2);
 \node[above right,fill=white,inner sep=1pt] at (0.15,5.15) {$p$};
 \node[right,fill=white,inner sep=1pt] at (\cc+0.1,0.0) {\small $\theta=0$};
 \node[left,fill=white,inner sep=1pt] at (-\cc+3.14159-0.1,0.0) {\small $\theta=\pi$};
 \node[fill=white,inner sep=1pt,font=\footnotesize] at (-3.4,-8.5) {$Ct=0.5$};
 \draw[->,thin] (-2.3,-8.0)--(-0.1,-5.25);
 \node[fill=white,inner sep=1pt,font=\footnotesize] at (3.4,-8.5) {$Ct=1.5$};
 \draw[->,thin] (2.3,-8.0)--(0.1,-6.0);
\end{tikzpicture}
\caption{Level curves $\Gamma_c$ of $f$ (thin dashed) and the source point $p=(0,5)$; the curve $\Gamma_{f(p)}$ through $p$ is thick. Open-headed double arrows: long axes of the indicatrix $E_x$, tangent to $\Gamma_c$; solid single arrows: $\nabla f$. Shaded bands: $\{|f-f(p)|\le Ct\}$ for $Ct=0.5$ (dark) and $Ct=1.5$ (light), schematic values. The disc $r<2$ is where $f$ is modified; level curves are not shown there.}
\label{fig:spirals}
\end{figure}

\FloatBarrier
\section*{Acknowledgments}
AI tools were used during the preparation of this manuscript to assist with proofreading and checking mathematical calculations. The graphics were created using AI tools. The final manuscript was fully reviewed, revised, and verified by the author, who retains sole responsibility for the final content and mathematical accuracy of this work.


\begin{thebibliography}{99}

\bibitem{AH} T.~Ariyoshi, M.~Hino, Small-time asymptotic estimates in local Dirichlet spaces, \emph{Electron. J. Probab.} \textbf{10} (2005), 1236--1259.

\bibitem{Ber} Yu.M.~Berezanskii, \emph{Selfadjoint operators in spaces of functions of infinitely many variables}, Translations of Mathematical Monographs, American Mathematical Society, Providence, RI, 1986.

\bibitem{CGT} J.~Cheeger, M.~Gromov, M.~Taylor, Finite propagation speed, kernel estimates for functions of the Laplace operator, and the geometry of complete Riemannian manifolds, \emph{J. Differential Geom.} \textbf{17}, no.~1 (1982), 15--53.

\bibitem{C} P.R.~Chernoff, Essential self-adjointness of powers of generators of hyperbolic equations, \emph{J. Funct. Anal.} \textbf{12}, no.~4 (1973), 401--414.

\bibitem{C2} P.R.~Chernoff, Schr\"odinger and Dirac operators with singular potentials and hyperbolic equations, \emph{Pacific J. Math.} \textbf{72}, no.~2 (1977), 361--382.

\bibitem{Ch} A.A.~Chumak, Self-adjointness of the Beltrami--Laplace operator on a complete paracompact Riemannian manifold without boundary, \emph{Ukrainian Math. J.} \textbf{25}, no.~6 (1974), 649--655.

\bibitem{CL} E.A.~Coddington, N.~Levinson, \emph{Theory of Ordinary Differential Equations}, McGraw-Hill, New York, 1955.

\bibitem{CS} T.~Coulhon, A.~Sikora, Gaussian heat kernel upper bounds via the Phragm\'en--Lindel\"of theorem, \emph{Proc. London Math. Soc.} \textbf{96}, no.~2 (2008), 507--544.

\bibitem{Davies} E.B.~Davies, \emph{Heat Kernels and Spectral Theory}, Cambridge Tracts in Mathematics, vol.~92, Cambridge University Press, Cambridge, 1989.

\bibitem{Folland} G.B.~Folland, \emph{Real Analysis: Modern Techniques and Their Applications}, 2nd ed., Wiley, New York, 1999.

\bibitem{FOT} M.~Fukushima, Y.~Oshima, M.~Takeda, \emph{Dirichlet Forms and Symmetric Markov Processes}, 2nd revised and extended ed., De Gruyter, Berlin, 2011.

\bibitem{Kato} T.~Kato, A remark to the preceding paper by Chernoff, \emph{J. Funct. Anal.} \textbf{12}, no.~4 (1973), 415--417.

\bibitem{NTX} S.-M.~Ngai, W.~Tang, Y.~Xie, Wave propagation speed on fractals, \emph{J. Fourier Anal. Appl.} \textbf{26} (2020), article~31.

\bibitem{O2} Yu.B.~Orochko, Finite propagation velocity and essential self-adjointness of certain differential operators, \emph{Funct. Anal. Appl.} \textbf{13}, no.~3 (1979), 235--237.

\bibitem{O} Yu.B.~Orochko, The property of global finite rate of propagation of a second-order elliptic differential expression, \emph{Differ. Uravn.} \textbf{18}, no.~10 (1982), 1764--1772 (in Russian).

\bibitem{P} M.A.~Perelmuter, Hyperbolic equation method in $L^p$-spaces, \emph{Ukrainian Math. J.} \textbf{41}, no.~1 (1989), 102--107.

\bibitem{PS} M.A.~Perelmuter, Yu.A.~Semenov, Finiteness of the rate of propagation of perturbations for hyperbolic equations, \emph{Ukrainian Math. J.} \textbf{36}, no.~1 (1984), 56--63.

\bibitem{Pov}  A.Ya.~Povzner, On the expansion of arbitrary functions in characteristic functions of the operator $-\Delta u+cu$, \emph{Mat. Sb. (N.S.)} \textbf{32(74)}, no.~1 (1953), 109--156 (in Russian).

\bibitem{RT} J.~Rauch, M.~Taylor, Essential self-adjointness of powers of generators of hyperbolic mixed problems, \emph{J. Funct. Anal.} \textbf{12}, no.~4 (1973), 491--493.

\bibitem{RB} F.S.~Rofe-Beketov, Necessary and sufficient conditions for a finite propagation rate for elliptic operators, \emph{Ukrainian Math. J.} \textbf{37}, no.~5 (1985), 547--549.

\bibitem{Sik} A.~Sikora, Riesz transform, Gaussian bounds and the method of wave equation, \emph{Math. Z.} \textbf{247}, no.~3 (2004), 643--662.

\bibitem{SW} E.M.~Stein, G.~Weiss, \emph{Introduction to Fourier Analysis on Euclidean Spaces}, Princeton Mathematical Series, no.~32, Princeton University Press, Princeton, N.J., 1971.

\end{thebibliography}
\end{document}